\documentclass[preprint,12pt]{elsarticle}
\usepackage[centertags]{amsmath}
\usepackage[utf8]{inputenc} 

\usepackage{hyperref}

\usepackage{amsmath}
\usepackage{amsfonts}
\usepackage{amssymb}
\usepackage{amssymb}
\usepackage{graphicx}
\usepackage{amssymb}
\usepackage{newlfont}
\usepackage{fancyhdr}
\usepackage{t1enc}
\usepackage{amsfonts}

\usepackage{ifpdf}

\newtheorem{thm}{Theorem}[section]

\newtheorem{cor}[thm]{Corollary}
\newtheorem{lem}[thm]{Lemma}
\newtheorem{prop}[thm]{Proposition}
\newtheorem{pr}[thm]{Properties}
\newproof{pf}{Proof}
\newtheorem{defn}[thm]{Definition}

\let \dis =\displaystyle

\usepackage{relsize}
\usepackage{amssymb}

\makeatletter
\def\ps@pprintTitle{%
	\let\@oddhead\@empty
	\let\@evenhead\@empty
	\def\@oddfoot{\reset@font\hfil\thepage\hfil}
	\let\@evenfoot\@oddfoot
} \makeatother

\usepackage{amssymb}
\usepackage{newlfont}
\usepackage{fancyhdr}
\usepackage{t1enc}
\usepackage{amsfonts}
\usepackage{cancel}
\usepackage{float}
\usepackage{caption}

\begin{document}

\begin{frontmatter}

\title{Some properties of q-Gaussian distributions}
	
\author[]{Ben Salah Nahla$^{a}$  \fnref{}}

\ead{bensalahnahla@yahoo.fr}
\cortext[cor1]{$^{a}$Corresponding author}
\address[els]{Laboratory of Probability and  Statistics, Faculty
of Sciences Sfax, University of Sfax. Tunisia}
\begin{abstract}
The q-Gaussian is a probability distribution generalizing the Gaussian one. In spite of a q-normal distribution is popular, there is a problem when calculating an expectation value with a corresponding normalized distribution and not a q-normal distribution itself.
  In this paper, two q-moments types called normalized and unormalized q-moments are introduced in details. Some properties of q-moments are given, and several relationships between them are established,
and some results related to q-moments are also obtained. Moreover, we show  that these new q-moments   may be regarded as a generelazation of the classical case for q = 1.
Firstly, we  determine the q-moments of q-Gaussian distribution. Especially, we give explicitely the kurtosis parameters. Secondly, we compute the expression of the q-Laplace transform of the q-Gaussian distribution. Finally, we study the distribution of sum of q-independent Gaussian distributions.
\end{abstract}

\end{frontmatter}


\section{Introduction and Preliminaries}\label{intro}
Several q-analogues of certain probability distributions have been recently investigated by many authors (\cite{A,1,2,3}). The q-distributions have been introduced in statistical physics for the characterization of chaos and multifractals. These distributions $f_{q}$ are a simple one parameter transformation of an original density function $f$ according to
$$f_{q}(x)=\dis\frac{f(x)^{q}}{\dis\int f(x)^{q} dx}$$
The parameter $q$ behaves as a microscope for exploring different regions of the measure $p$: for $q > 1$, the more singular regions are amplified, while for $q < 1$ the less singular regions are accentuated. 
\\Tsallis distributions (q-distributions) have encountered a large success because of their remarkable agreement with experimental data, see \cite{A1, A2, A3,I, II}, and references therein.
In particular, the q-Gaussian distribution is also well-known as a gneralisation of the Gaussian, or the Normal distribution. This distribution can also represent the heavy tailed distribution such as the Student-distribution or  the distribution with bounded support such as the semicircle  of Wigner. For these reasons, the q-Gaussian distribution has been applied  in the fields of statistical mechanics, geology, finance, and machine learning.  Admitting the q-normal distribution is in demand as above, there exists a problem to calculate the expectation value with a corresponding q-distribution not a q-normal distribution itself. 
But we have an amazing property such that an escort distribution obtained by a q-normal distribution with a parameter $q$ and a variance is another q-normal distribution with a different value of $q$ and a scaled variance. Then calculating an expectation value with an escort distribution corresponds to calculating the expectation value with another associated q-normal distribution, but it gets even the question why an expectation value should be calculated by another q-normal distribution. We call the procedure to get another q-normal distribution from a given q-normal distribution through an escort distribution   proportion.

Furthermore, we target attention on q-Gaussians, an essential tool of q-statistics \cite{26}, that was not discussed in \cite{25}. The q-Gaussian behavior is often detected in quite distinct settings \cite{26}.

It is well known that, in the literature, there are two types of q-Laplace transforms, and they are studied in detail by several authors (\cite{A4,A5}, etc.).
Recently Tsallis et al. have been interested in calculating  the Fourier transform of q-Gaussian   and  have proved a generalization of the central limit theorem for $1 \leq q < 3$. The case $q < 1$ requires essentially different technique, therefore we leave it for a separate paper. In this paper,  we propose new definitions of the q-laplace transform of some probability distributions. These results are motivated by recent developments in the calculation of Fourier transforms, where new formulas have been defined \cite{3}.

In this article, we develop our results into four sections. In Section 2, we recall some known definitions and notations from the q-theory.

In Section 3, we give definitions of some q-analogues of  mean  and variance. In Section 4, we introduce the q-Gaussian distribution  includes some properties.  
In Section 5, we give the news formula of Laplace transform and we treat kurtosis both in its standard definition and in q statistics, namely q-kurtosis. In Section 6, we estimate the q-mean and q-variance.
\\We start with definitions and facts from the q-calculus.
\section{q-theory calculus}
Assume that $q$  be a fixed number satisfying $q \in [0,1].$ If is a classical object,say, its q-version is defined by $[x]_{q}=\dis \frac{1-q^{n}}{1-q}.$ As is well know, the q-exponential and the q-logarithm, which are denoted by $e_{q}(x)$ and $\ln_{q}(x),$ are respectively defined as $e_{q}(x)=[1+(1-q)x]_{+}^{\frac{1}{1-q}}$ and $\ln_{q}(x)=\frac{x^{1-q}-1}{1-q},(x>0).$
For q-exponential, the relations $e_{q}^{x \otimes_{q} y}=e_{q}^{x}e_{q}^{y}$ and $e_{q}^{x +y}=e_{q}^{x} \otimes_{q}  e_{q}^{y}$ hold true. These relations can be rewritten equivalently as follows:
$\ln_{q} (x \otimes_{q} y)=\ln_{q} (x)+\ln_{q} (y),$ and 
$\ln_{q} (x y)=\ln_{q} (x) \otimes_{q} \ln_{q} (y).$ 

\vskip0.5cm A q-algebra can also be defined in \cite{A} by applying the generalized operation
for sum and product:
$$x \oplus_{q} y = x + y + (1 - q)xy ,$$ 
$$x \otimes_{q} y=[x^{1-q}+y^{1-q}-1]_{+}^{1-q},$$
 with the following neutral and inverse elements:
 $$x \otimes_{q} (x)_{q}=0,  \hbox{ with } (x_{q}) =x [ 1+(1-q)x]^{-1}$$

$$x \otimes_{q}(x^{-1})_{q} =1, \hbox{ with }:(x^{-1})_{q}  =[2-x^{1-q}]_{+}^{\frac{1}{1-q}} .$$

For the new algebraic operation, q-exponential and q-logarithm have the following properties:
\begin{pr}

\begin{enumerate}
\item $e_{q}^{x} e_{q}^{y}=e_{q}^{x \otimes_{q} y}$

\item $e_{q}^{x} \otimes_{q} e_{q}^{y}=e_{q}^{x+y} $

\item $\log_{q}(xy)=\log_{q}(x) \otimes_{q} \log_{q}(y)$

\item $\log_{q}(x \otimes_{q} y)=\log_{q}(x) +\log_{q}(y)$
\end{enumerate}

\end{pr}

It can be easily proved that the operation $\otimes_{q}$ and $\oplus_{q}$ satisfy commutativity and associativity. For the operator $\oplus_{q},$ the identity additive is 0, while for the operator $\otimes_{q}$ the identity multiplicative is 1 \cite{C}. Two distinct mathematical tools appears in the study of physical phenomena in the complex media which is characterized by singularities in a compact space.

From the associativity of $\oplus_{q}$ and $\otimes_{q}$, we have the following formula :
$$t \oplus_{q} t \oplus_{q}....\oplus_{q} t=\frac{1}{1-q}\{[1+(1-q)t]^{n}-1\}$$
$$t^{\otimes_{q} n}=t \otimes_{q} t \otimes_{q} ...\otimes_{q} t={n  t^{1-q}-(n-1)}^{\frac{1}{1-q}}.$$
The real space vector with regular sum and product operations $\mathbb{R}(+, \times)$ is a field, and the $\mathbb{R}(\oplus_{q}, \otimes_{q})$ defines a quasi-field.
\section{q-mean and q-variance values }
Let $q$ be a real number and $f$ be  a properly normalized probability density with $supp f \subseteq \mathbb{R}$ of some random variable $X$ such that the quantity
$$\dis\int_{-\infty}^{+\infty} f(x)dx = 1.$$
The mean $m$ is defined, of a given $X$, as follows
$$E(X)=m=\dis \int_{-\infty}^{+\infty} x f(x)dx.$$
The variance $V$ is defined, of a given $X$, as follows
$$V(X)=\dis\int_{-\infty}^{\infty}  (x-m)^{2} f(x)dx.$$
The unnormalized q-moments , of a given $X$,  is defined as 
$$E_{q}(X)=m_{q}=\dis \int_{-\infty}^{+\infty} x [f(x)]^{q}dx.$$
Similarly, the  unnormalized q-variance, $\sigma^{2}_{2q-1}$ is defined analogously to the usual second order central moment, as
$$V_{2q-1}(X)=\sigma^{2}_{2q-1}= \dis \int_{-\infty}^{+\infty} (x-m_{2q-1})^{2} [f(x)]^{2q-1}dx.$$
On the other hand, we denote by $f_{q}(x)$ the normalized density (see e.g. \cite{A6}) and defined as
$$f_{q}(x)=\dis\frac{[f(x)]^{q}}{\nu_{q}(f)}.$$
where
$$\nu_{q}(f)=\dis\int_{-\infty}^{+\infty} [f(x)]^{q} dx<\infty.$$
The normalized q-mean values, of a given $X$, is

$$\overline{E_{q}}(X)=\overline{m}_{q}=\dis \int_{-\infty}^{+\infty} x f_{q}(x)dx$$
The normalized $q$-variance values, of a given $X$, is
$$\overline{V}_{2q-1}(X)=\overline{\sigma}^{2}_{2q-1}=\dis \int_{-\infty}^{+\infty} (x-m_{2q-1})^{2} f_{2q-1}(x)dx.$$

\section{q-Gaussian  distribution}
In  this section, we review the q-Gaussian
distribution, or the q-normal distribution according to Furuichi \cite{A}  and Suyari \cite{B}. Let $\beta$ be a positive number. We call the q-Gaussian $N_{q}(m,\sigma^{2})$ with parameters $m$ and $\sigma^{2}>0$ if its density function is defined by
$$f(x)=\dis\frac{\sqrt{\beta}}{ \sigma C_{q}} e_{q}^{ -\dis\frac{\beta {(x-m)}^{2}}{\sigma^{2}}}, x \in \mathbb{R}$$
with  $q<3, q \neq 1$; and $C_{q}$ is the normalizing constant, namely
$$C_{q}=\int_{-\infty}^{\infty} e_{q}^{-x^{2}}dx=\left\{%
\begin{array}{ll}
& \Big{(}\dis\frac{1}{q-1}   \Big{)} ^{\frac{1}{2}}  B\Big{(}  \dis\frac{3-q}{2(q-1)},\frac{1}{2}   \Big{)} 
\hskip2cm     1  \textless q \textless 3                 \\
&\sqrt{\pi}, \hskip5cm    q=1\\
&\Big{(}\dis\frac{1}{1-q}   \Big{)} ^{\frac{1}{2}}  B\Big{(}  \dis\frac{2-q}{1-q},\frac{1}{2}   \Big{)}   \hskip3cm     -\infty \textless q \textless 1	\\
\end{array}%
\right.$$
and $B(a,b)$ denotes the beta function.
The widht parameter of the distribution is characterized by
$$\beta=\frac{1}{3-q}.$$
Denote a general q-Gaussian random variable $X$ with parameters  $m$ and $ \sigma^{2}$ as $X \sim N_{q}(m,\sigma^{2}),$ and call the special case of $m=0$ and  $\sigma^{2}=1$ a standard q-Gaussian $Y \sim N_{q}(0,1)$. The density of the satndard q-Gaussian distribution may then be written as
$$N_{q}(0,1)(y)=\dis\frac{\sqrt{\beta}}{C_{q}} e_{q}^{-\beta y^{2}}.$$ 
\\Note that, if	
\begin{equation} \label{5} Y \sim N_{q}(0,1) \hbox{ then } X=m+ \sigma Y \sim N_{q}(m,\sigma^{2}).\end{equation}
If we change the value of $q$, we can represent various types of distributions. The $q$-Gaussian distribution represents the usual Gaussian distribution when $q = 1,$ has compact support for $q < 1,$ and turns asymptotically as a power law for $1 \leq q<3.$ For $3 \leq q,$ the form given  is not normalizable. The usual variance (second order moment) is finite for $q < \frac{5}{3}$ , and, for the standard q-Gaussian $N_{q}(0,1)$, is given by $V(Y) =\dis\frac{3-q}{5-3q}.$ The usual variance of the q-Gaussian diverges for $\frac{5}{3} \leq q < 3$ , however the q-variance remains finite for the full range $\infty < q < 3 $, equal to unity for the standard q-Gaussian.
Finally, we can easily check that there are relationships
between different values of q. For example, \begin{equation}\label{zz} e_{q}^{-y^{2}}=(e_{2-\frac{1}{q}}^{-q y^{2}})^{\frac{1}{q}}.\end{equation}
In this section, we consider the q-analogues of the Laplace transform, which we call the q-Laplace transform, and investigate some of its properties.
\section{New q-Laplace transforms}
From now, we assume that $1  \leq q < 3$. For these values of $q$  we introduce the q-Laplace transform $L_{q}$ as an operator, which coincides with the Laplace transform if $q = 1$. Note that the q-Laplace transform is defined on the basis of the q-product and the q-exponential, and, in contrast to the usual Laplace transform, is a nonlinear transform for $q \in  (1, 3).$

The q-Laplace transform of a random variable $X$ with density function $f$ is defined by the formula 

$$L_{q}(X)(\theta)=\dis\int_{supp f} e_{q}^{ \theta x} \otimes_{q} f(x)dx,.$$
where the integral is understood in the Lebesgue sense. 

The following lemma establishes the expression of the q-Laplace transform in terms of the standard product, instead of the q-product.

\begin{lem} \label{n}
	The q-Laplace transform of a random variable $X$ with density $f$ is expressed as 
	\begin{equation} \label{55}  L_{q}(X)(\theta)=\dis\int_{-\infty}^{\infty}  f(x) e_{q}^{ \theta x (f(x))^{q-1}  }dx.\end{equation} 
	
\end{lem}

\begin{pf}
	For $x \in supp f,$ we have
	\begin{equation} \label{66}  e_{q}^{i\theta x } \otimes_{q} f(x)=[1+(1-q)\theta x+(f(x))^{1-q}-1]^{\frac{1}{1-q}}
	=f(x)[1+(1-q)\theta x (f(x))^{q-1}]^{\frac{1}{1-q}} \end{equation} 
	
\end{pf}
Integrating both sides of Eq. ($\ref{66}$)  we obtain  $(\ref{55}).$  
\\Let $X$ be a random variable defined on the probability space $(\Omega,F, P)$ with density function $f \in L^{q}$.
It can be verified that the derivatives of the q-Laplace transform $L_{q}(X)(\theta)$ are closely related to an appropriate set of unnormalized q-moments of the original probability density. Assume that $L_{q}(X)(\theta)< +\infty$ in a neighbor of 0.

Indeed, the first few low-order derivatives (including the zeroth order) are given by
$$L_{q}(X)(0)=1$$
$$\frac{\partial L_{q}(X)(\theta)}{\partial \theta}_{\Big{|} \theta=0}= \dis\int_{-\infty}^{\infty}  x (f(x))^{q}  dx=E_{q}(X)$$
$$\frac{\partial^{2} L_{q}(X)(\theta^{2})}{\partial \theta}_{\Big{|} \theta=0}=q \dis\int_{-\infty}^{\infty}  x^{2} (f(x))^{2q-1}  dx=q E_{2q-1}(X^{2}) $$
$$\frac{\partial^{3} L_{q}(X)(\theta^{3})}{\partial \theta}_{\Big{|} \theta=0}=q(2q-1) \dis\int_{-\infty}^{\infty}  x^{3} (f(x))^{3q-2}  dx=q(2q-1) E_{3q-2}(X^{3}).$$
$$\frac{\partial^{4} L_{q}(X)(\theta^{4})}{\partial \theta}_{\Big{|} \theta=0}=q(2q-1)(3q-2) \dis\int_{-\infty}^{\infty}  x^{4} (f(x))^{4q-3}  dx=q(2q-1)(3q-2) E_{4q-3}(X^{4}).$$
The general n-derivative is
$$\frac{\partial^{n} L_{q}(X)(\theta^{n})}{\partial \theta}_{\Big{|} \theta=0}=\prod_{m=0}^{n-1} (1+m(q-1)) \dis\int_{-\infty}^{\infty}  x^{n} (f(x))^{1+n(q-1)}  dx, n=1,2,3....$$

Note that, in the case $n=1$ the first derivative of the Laplace transform corresponds to $E_{q}(X).$

\begin{prop}\label{4}
	Let $1 \leqslant q < 3$ and let $X$ be a random variable following  a q-Gaussian distribution $N_{q}(m,\sigma^{2} )$
 then $E(X)=m$ and $V(X)=\dis\frac{3-q}{5-3q} \sigma^{2}$ with $1\leq q <\frac{5}{3}.$
\end{prop}
\begin{pf}
\begin{enumerate}	
	\item 	The first moment, of a given $X,$ is
\begin{align*}
	&\begin{aligned}{}
E(X)&=\dis\frac{\sqrt{\beta}}{ \sigma C_{q}} \dis\int_{-\infty}^{\infty} x  e_{q}^{ -\dis\frac{\beta {(x-m)}^{2}}{\sigma^{2}}}dx\\
    &= \dis\frac{\sqrt{\beta}}{  C_{q}}  \dis\int_{-\infty}^{\infty} (\sigma y +m)  e_{q}^{ -\beta {y}^{2}}dy   \\
    &=	\sigma \dis\frac{\sqrt{\beta}}{  C_{q}} \dis \int_{-\infty}^{\infty}  y  e_{q}^{ -\beta {y}^{2}}dy +m \dis\frac{\sqrt{\beta}}{ C_{q}}	\dis \int_{-\infty}^{\infty}   e_{q}^{ -\beta {y}^{2}}dy =m \\
\end{aligned}&&
\end{align*} 
\item The second order moment  of the standard Gaussian $N_{q}(0,1)$ is computed as
	\begin{align*}
	&\begin{aligned}{}
	E(Y^{2})&=\dis\frac{\sqrt{\beta}}{ C_{q}} \int_{-\infty}^{\infty} y^{2}            e_{q}^{-\beta y^{2}} dy\\        
	&=\dis\frac{1}{ 2\sqrt{\beta} (2-q) C_{q}}    \dis\int_{-\infty}^{\infty} (e_{q}^{-\beta y^{2}})^{2-q}  dy\\
	&=\dis\frac{1}{ 2\sqrt{\beta} (2-q) C_{q}}    \dis\int_{-\infty}^{\infty} e_{q_{1}}^{-\beta (2-q)y^{2}} dy, q_{1}=\frac{1}{2-q}\\
	& \hbox{ The substitution }  \beta_{1} z^{2} =\beta (2-q)y^{2}, \beta_{1}=\frac{1}{3-q_{1}} \\
	&=\dis\frac{\sqrt{\beta_{1}} }{ 2 \beta (2-q)^{\frac{3}{2}} C_{q}}    \dis\int_{-\infty}^{\infty} e_{q_{1}}^{-\beta_{1} z^{2}} dz\\
	&=\dis\frac{C_{q_{1}}}{ 2 \beta (2-q)^{\frac{3}{2}} C_{q}}    \dis\int_{-\infty}^{\infty} \frac{\sqrt{\beta_{1}}}{ C_{{q}_{1}}} e_{q_{1}}^{-\beta_{1} z^{2}} dz\\
	&=\dis\frac{1 }{ 2 \beta_{1} (2-q)^{\frac{3}{2}}} \frac{C_{q_{1}}}{C_{q}}\\                	
	\end{aligned}&&
	\end{align*} 
The condition $1 \leqslant q<3$ implies that $1 \leqslant q_{1}<\frac{5}{3}.$
By using the identity  $B(x+1,y)=\frac{x}{x+y} B(x,y)$ we obtain  the ratio between $C_{q}$ and $C_{q_{1}}$  as
\begin{equation} \label{2}	\dis\frac{C_{q_{1}}}{C_{q}}=\dis\frac{2 (2-q)^{\dis\frac{3}{2}}}{5-3q} \end{equation} 
\end{enumerate} 
By applying the formula $V(X)=E(X^{2})-(E(X))^{2}$, we obtain the result	
\end{pf}

In this theorem, we give  the average of the fourth power of the standardized deviations from the q-mean.

In this theorem, we determine the q-kurtosis of q- Gaussian distribution ($N_{q}(0,1)$).

\begin{thm} \label{555}
	Let $1 \leqslant q < 3$ and let $X$ be a random variable following  a q-central Gaussian distribution $N_{q}(0,1)$, then the coefficient of kurtosis is 
	$$Kurt[Y]=\dis\frac{E(Y^{4})}{(E(Y^{2}))^{2}}=\dis\frac{3 (5-3q)}{(7-5q)}, 1\leq q<\frac{7}{5}$$
\end{thm}
	
\begin{pf}
	\begin{enumerate}
		\item
	The fourth central moment  moment of the standard Gaussian $N_{q}(0,1)$ is computed as
		
	\begin{align*}
	&\begin{aligned}{}
	E(Y^{4})&=\dis\frac{\sqrt{\beta}}{ C_{q}} \int_{-\infty}^{\infty} y^{4}                     e_{q}^{-\beta y^{2}} dy\\
	&=\dis\frac{3}{ 2\sqrt{\beta} (2-q) C_{q}(1)}    \dis\int_{-\infty}^{\infty}  y^{2} (e_{q}^{-\beta y^{2}})^{2-q}  dy\\
	&\hbox{ The substitution }  \beta_{1} z^{2} =\beta (2-q)y^{2}\\
	&=\dis\frac{3 \beta_{1} }{ 2 \beta^{2} (2-q)^{\frac{5}{2}}} \frac{C_{q}}{C_{q}}  \dis\int_{-\infty}^{\infty}  \frac{\sqrt{\beta_{1}}}{ C_{q_{1}}} z^{2} e_{q_{1}}^{-\beta_{1} z^{2}} dz \hbox{ with } 1 \leq q_{1}=\frac{1}{2-q}<\frac{5}{3}\\
	&=\dis\frac{3 \beta_{1} }{ 2 \beta^{2} (2-q)^{\frac{5}{2}}} \frac{C_{q_{1}}}{C_{q}}  E_{q_{1}}(Z^{2})\\
	&\hbox{ Using equation  \ref{2}, we obtain}  \\
	&=\dis\frac{3(3-q)^{2}}{(5-3q)(7-5q)} ,1 \leq q<\frac{7}{5} \\
	\end{aligned}&&
	\end{align*} 
	According to Proposition  $\ref{4},$ we obtain the result.

\end{enumerate}

\end{pf}

For $1 \leq q<\frac{7}{5}$ a value greater than $\frac{3(3-q)^{2}}{(5-3q)(7-5q)}$ indicates a leptokurtic distribution; a values less than $\frac{3(3-q)^{2}}{(5-3q)(7-5q)}$ indicates a platykurtic distribution. For the sample estimate $X$, $\frac{3(3-q)^{2}}{(5-3q)(7-5q)}$ is subtracted so that a positive value indicates leptokurtosis and a negative value indicates platykurtosis. 
\begin{thm}
Let $1 \leqslant q < 3$ and let $X$ be a random variable following  a q-Gaussian distribution $N_{q}(m,\sigma^{2})$, then 
		\begin{enumerate}	
			\item $L_{q}(X)(\theta)=\Big{(} e_{q}^{\theta m{a^{q-1}  -\dis\frac{\theta^{2} a^{2q-2} \sigma^{2} }{4 \beta}}}\Big{)}^{\frac{3-q}{2}},\hbox{ with } a=\dis\frac{\sqrt{\beta}}{ \sigma C_{q}} \hbox{  and }\theta \in \mathbb{R}$
		\item $E_{q}(X)=\dis\int_{-\infty}^{\infty}  x (f(x))^{q}  dx=\dis\frac{ m (3-q)^{\dis\frac{3-q}{2}}}{ 2(\sigma C_{q})^{q-1}}$
		\item $E_{2q-1}(X^{2})=\dis\int_{-\infty}^{\infty}  x^{2} (f(x))^{2q-1} dx =\dis\frac{1}{4 q (3-q)^{q-2}(\sigma C_{q})^{2q-2}}[(3-q)\sigma^{2}+ (q+1)m^{2} ]$
	\end{enumerate}	
	
\end{thm}

\begin{pf} 
\begin{enumerate}	
\item	From definition of q-Laplace transform, it following by denote  $a=\dis\frac{\sqrt{\beta}}{\sigma C_{q}}$ 
\begin{align*}
L_{q}(X)(\theta)&= \dis\int_{-\infty}^{\infty}    e_{q}^{\theta x} \otimes_{q} a  e_{q}^{\dis\frac{-\beta(x-m)^{2}}{\sigma^{2}}}dx\\
&\hbox{(by appliying lemma \ref{n})}\\	
&= a \dis\int_{-\infty}^{\infty}    e_{q}^{\dis\frac{-\beta (x-m)^{2}}{\sigma^{2}}} e_{q}^{ \theta x   { a^{q-1}  \Big{(}e_{q}^{ \dis\frac{-\beta (x-m)^{2}}{\sigma ^{2}}}\Big{)}^{q-1}}}dx\\
&= a \dis\int_{-\infty}^{\infty} e_{q}^{\theta x a^{q-1}} \otimes_{q} e_{q}^{\dis\frac{-\beta (x-m)^{2}}{\sigma^{2}}}    dx\\
&= a\dis\int_{-\infty}^{\infty} e_{q}^{\dis\frac{-\beta (x-m)^{2}}{\sigma^{2}}+\theta x a^{q-1}}  dx\\
&= a \dis\int_{-\infty}^{\infty} e_{q}^{{\dis\frac{-\beta (x-m)^{2}}{\sigma ^{2}}}+   \dis\frac{\beta}{\sigma^{2}} (\dis\frac{\sigma^{2} \theta a^{q-1}}{2 \beta})^{2} +\theta a^{q-1} m}  dx\\
&= a \dis\int_{-\infty}^{\infty} e_{q}^{\dis\frac{\beta}{\sigma^{2}} (\dis\frac{\sigma^{2} \theta a^{q-1}}{2 \beta})^{2} +\theta a^{q-1} m} e_{q}^{{\dis\frac{-\beta (x-m)^{2}}{\sigma^{2}}}  (e_{q}^{\dis\frac{\beta}{\sigma^{2}} (\dis\frac{\sigma^{2} \theta a^{q-1}}{2 \beta})^{2} +\theta a^{q-1} m})^{q-1}}  dx\\
&= a e_{q}^{\dis\frac{\beta}{\sigma^{2}} (\dis\frac{\sigma^{2} \theta a^{q-1}}{2 \beta})^{2} +\theta a^{q-1} m}  \dis\int_{-\infty}^{\infty} e_{q}^{{\dis\frac{-\beta \gamma (x-m)^{2}}{\sigma^{2}}}}  dx\\
&= \sigma_{1} a  \frac{C_{q}}{\sqrt{\beta}}e_{q}^{\dis\frac{\beta}{\sigma^{2}} (\dis\frac{\sigma^{2} \theta a^{q-1}}{2 \beta})^{2} +\theta a^{q-1} m} \\
&= \frac{e_{q}^{\dis\frac{\beta}{\sigma^{2}} (\dis\frac{\sigma^{2} \theta a^{q-1}}{2 \beta})^{2} +\theta a^{q-1} m}}{\sqrt{e_{2}^{(q-1)(\dis\frac{\beta}{\sigma^{2}} (\dis\frac{\sigma^{2} \theta a^{q-1}}{2 \beta})^{2} +\theta a^{q-1} m)}}} \\
\end{align*}           
where $\gamma= e_{q}^({\dis\frac{\beta}{\sigma^{2}} (\dis\frac{\sigma^{2} \theta a^{q-1}}{2 \beta})^{2} +\theta a^{q-1} m})^{q-1} \hbox{and} \sigma_{1}^{2}=\dis\frac{\sigma^{2}}{\gamma}.$       \\Hence, applying again Lemma  $\ref{n},$ we have
\begin{align*}   
L_{q}(X)(\theta)&=e_{q}^{ \Big{(}  \theta m{ a^{q-1} -\dis\frac{\theta^{2} a^{2q-2} \sigma^{2} }{4 \beta}\Big{)}}^{\frac{3-q}{2}}}    \\&=e_{q_{1}}^{{\dis\frac{ \theta m(3-q)}{2} a^{q-1}   -\dis\frac{\theta^{2} a^{2q-2} \sigma^{2}(3-q) }{8 \beta}}}    \\
\end{align*}
\item We compute the first and the second derivative of $L_{q}(X),$ with respect to $\theta$, we obtain 
$$E_{q}(X)=L'_{q}(X)(0)$$
$$ E_{2q-1}(X^{2})=\frac{1}{q}L"_{q}(0)$$
\end{enumerate}
\end{pf}


Intending to interpret these moments, we consider  a wonderful property such that an escort distribution obtained by a q-normal distribution with variance $\sigma^{2}$ is
equivalent to another q-normal distribution with $q_{1} = 2 -\frac{1}{q}$ and a variance $\dis\frac{3-q}{q+1} \sigma^{2}$ with $1 \leq q<3$.
\begin{prop}
 Let $X$ be a random variable following  a q-Gaussian distribution $N_{q}(m,\sigma^{2})$,then
\begin{enumerate} 
\item $\overline{E}_{q}(X)=\dis \int_{-\infty}^{+\infty} x \dis\frac{[f(x)]^{q}}{\nu_{q}(f)}dx=m$
\item $\overline{V}_{2q-1}(X)=\sigma^{2}_{2q-1}=\dis \int_{-\infty}^{+\infty} (x-m)^{2} \dis\frac{[f(x)]^{2q-1}}{\nu_{2q-1}(f)}dx=\dis\frac{3-q}{q+1} \sigma^{2}$
\end{enumerate}

\end{prop}

\begin{pf} 
\begin{enumerate} 	
\item Let's begin with observing a following proportion on a given q-Normal distribution:
$$ q_{3}=\frac{2q-1}{q}, \sigma_{3}^{2}=\frac{\beta_{3}}{(q) \beta} \sigma^{2}=\frac{3-q}{q+1} \sigma^{2}$$	
Under this relations, we have
\begin{align*} 
\dis\frac{\sqrt{\beta}}{ \sigma C_{q}}   \frac{\Big{(}e_{q}^{ -\dis\frac{\beta {(x-m)}^{2}}{\sigma^{2}}}\Big{)}^{q}}{\nu_{q}}                            &\propto e_{q_{3}}^{-\dis\frac{\beta {(x-m)}^{2}}{\sigma^{2}}\frac{q \beta}{\beta_{3}}}\\
                &\propto N_{q_{3}}(m,\sigma_{3}^{2})\\
\end{align*} 	
Therefore, 
$$\overline{E}_{q}(X)=\dis \int_{-\infty}^{+\infty} x \dis\frac{[f(x)]^{q}}{\nu_{q}(f)}dx=\int_{-\infty}^{+\infty} x N_{q_{3}}(m,\sigma_{3}^{2})(x)dx$$

By applying proposition $\ref{4},$ we obtain the result.

\item The escort function is proportional to the q-gaussian $N_{q_{4}}(m,\sigma^{2})$
$$ q_{4}=\frac{3q-2}{2q-1}, \sigma_{4}^{2}=\frac{\beta_{4}}{(2q-1) \beta} \sigma^{2}=\frac{3-q}{3q-1} \sigma^{2}.$$
Hence
\begin{equation}
\begin{split}
\dis\frac{\sqrt{\beta}}{ \sigma C_{q}}   \frac{\Big{(}e_{q}^{ -\dis\frac{\beta {(x-m)}^{2}}{\sigma^{2}}}\Big{)}^{q}}{\nu_{2q-1}}&\propto e_{q_{4}}^{-\frac{(x-m)^{2}}{\sigma^{2}}    \frac{(2q-1) \beta}{\beta_{4}}}\\
&\propto N_{q_{4}}(m,\sigma_{4}^{2})\\
\end{split}
\end{equation}	
Then
$$\overline{V}_{2q-1}(X)=\dis \int_{-\infty}^{+\infty} (x-m)^{2} \dis\frac{[f(x)]^{2q-1}}{\nu_{2q-1}(f)}dx=\dis \int_{-\infty}^{+\infty} (x-m)^{2} N_{q_{4}}(m, \sigma_{4}^{2})(x)dx=\frac{3-q_{4}}{5-3q_{4}} \sigma_{4}^{2}$$	
By applying proposition $\ref{4}$we obtain the result.
\end{enumerate}
\end{pf}

In rhis theorem, we prove that for $q<\frac{3}{5}$ a value greater than $\frac{3 (q+1)^{2}}{(5q-3)(3q-1)}$ indicates a leptokurtic distribution; a values less than $\frac{3 (q+1)^{2}}{(5q-3)(3q-1)}$ indicates a platykurtic distribution. For the sample estimate $X$, $\frac{3 (q+1)^{2}}{(5q-3)(3q-1)}$ is subtracted so that a positive value indicates leptokurtosis and a negative value indicates platykurtosis. 
\begin{thm}

Let $Y$ be a random variable following  a q-Central Gaussian distribution $N_{q}(0,1)$ then  the coefficient of normalized kurtosis  is 
$$\overline{Kurt}[Y]=\dis\frac{\overline{E}(Y^{4})}{(\overline{E}(Y^{2}))^{2}}=\dis\frac{3 (q+1)^{2}}{(5q-3)(3q-1)},       1 \leq q<\frac{3}{5}$$
\end{thm}

\begin{pf}
	
Checking a following proportion on a given q-Normal distribution:
$$q_{1}=\frac{5q-4
}{4q-3}, \sigma_{1}^{2}=\frac{\beta_{1}}{(4q-3) \beta}.$$
\begin{align*}
(\dis\frac{\sqrt{\beta}}{ \sigma C_{q}})^{4q-3} e_{q}^{ -\dis\frac{\beta {x}^{2}}{\nu_{4q-3}}}&\propto e_{q_{1}}^{-x^{2}\beta_{1} \frac{(4q-3) \beta}{\beta_{1}}   }\\
&\propto N_{q_{1}}(0, \sigma_{1}^{2}) \\
\end{align*}
Then
$\overline{E}(Z^{4})=\dis\int_{-\infty}^{+\infty} z^{4} (\dis\frac{\sqrt{\beta}}{ \sigma C_{q}})^{4q-3} e_{q}^{ -\dis\frac{\beta {z}^{2}}{\nu_{4q-3}}}=\int_{-\infty}^{+\infty} z^{4}N_{q_{1}}(0,\sigma_{1}^{2})(z)dz.$	
\\Hence, according to Theorem $\ref{555}$, we obtain $\overline{E}(Z^{4})=\dis\frac{3(3-q_{1})^{2}}{(5-3q_{1})(7-5q_{1})} \sigma_{1}^{4}; Z \sim N_{q_{1}}(0, \sigma_{1}^{2}).$
 we obtain
$$\overline{E}(Z^{4})=\frac{3(3-q)^{2}}{(5q-3)(3q-1)}.$$

Furthermore by observing a following proportion on a given q-Normal distribution:
$$q_{2}=\frac{3q-2
}{2q-1}, \sigma_{2}^{2}=\frac{\beta_{2}}{(2q-1) \beta}$$

\begin{equation}
\begin{aligned}
(\dis\frac{\sqrt{\beta}}{ \sigma C_{q}})^{2q-1} e_{q}^{ -\dis\frac{\beta {x}^{2}}{\nu_{2q-1}}}
&\propto& e_{q_{2}}^{-x^{2}\beta_{1} \frac{(2q-1) \beta}{\beta_{2}}   }&\\
&\propto &N_{q_{2}}(0, \sigma_{2}^{2})&\\
\end{aligned}
\end{equation}

Then
$$\overline{E}(Z^{2})=\dis\int_{-\infty}^{+\infty} z^{2} (\dis\frac{\sqrt{\beta}}{ \sigma C_{q}})^{2q-1} e_{q}^{ -\dis\frac{\beta {z}^{2}}{\nu_{2q-1}}}=\int_{-\infty}^{+\infty} z^{2}N_{q_{2}}(0,\sigma_{1}^{2})(z)dz.$$	
Hence, according to proposition $\ref{4}$, we obtain
$$\overline{E}(Z^{2})=\frac{(3-q_{2}) \sigma_{2}^{2}}{(5-3q_{2})}; Z \sim N_{q_{2}}(0, \sigma_{2}^{2})$$
 
$$\overline{E}(Z^{2})=\frac{3-q}{q+1}.$$
\end{pf}
q-estimator for random variables are arising from non-extensive statistical mechanics. In this section, we will estimate the q-mean and q-variance using the notions of q-Laplace transform, q-independence.
\section{Estimator of q-Mean and q-variance}
	\begin{defn} 
	Two random variables $X_{1}$ and $X_{2}$ are said to be q-independent if 
	$$L_{q}(X_{1}+X_{2})(\theta)=L_{q}(X_{1})(\theta) \otimes_{q} L_{q}(X_{2})(\theta).$$
\end{defn}

\begin{defn}	
Let $X_{n}$ be a sequence of identically distributed random variables and $m= E(X_{1}).$
Denote $S_{n}= \dis\sum_{k=1}^{n} X_{k}..$ By definition 
 $X_{k}, k=1,2,3,...,$ is said to be q-independent of the first type (or q-i.i.d.) if for all $n=2,3,4,..., $ the relations $$L_{q}[S_{n} -n m](\theta)= L_{q}[X_{1} -m](\theta)\otimes_{q} ....\otimes_{q} L_{q}[X_{n}-m](\theta)$$ hold.	
\end{defn}

\begin{prop}
Let  $X_{1} $ and $X_{2}$ be tow q-independent  random variables following respectively   $N_{q}(m_{1},\sigma_{1}^{2})$ and  $N_{q}(m_{2},\sigma_{2}^{2})$.
Then $$X_{1}+X_{2} \curvearrowright N_{q}(m_{1}+m_{2}, \sigma_{1}^{2}+\sigma_{2}^{2})$$
\end{prop}

\begin{pf}
	\begin{align*}
L_{q}(X_{1}+X_{2})(\theta)&=L_{q}(X_{1})(\theta) \otimes_{q} L_{q}(X_{2})(\theta)\\
                             &=e_{q_{1}}^{{  \dis\frac{3-q}{2} (\sigma_{1}^{2}+\sigma_{2}^{2})  \dis\frac{\theta^{2} a^{2q-2}}{4 \beta}}+        \dis\frac{3-q}{2} (\sigma_{1}^{2}+\sigma_{2}^{2}) \theta a^{q-1} }   \\
                           &=\Big{(}e_{q}^{{ (\sigma_{1}^{2}+\sigma_{2}^{2})  \dis\frac{\theta^{2} a^{2q-2}}{4 \beta}}+         (\sigma_{1}^{2}+\sigma_{2}^{2}) \theta a^{q-1} }\Big{)}^{\dis\frac{3-q}{2}}   \\
                           \end{align*}
\end{pf}
Note that if $X_{1}$ and $X_{2}$ are q-Gaussian  and q-independent random variables with distributions  $N_{q}(m_{1}, \sigma_{1}^{2} )$ and $N_{q}(m_{2}, \sigma_{2}^{2} )$ respectively then $$V(X_{1}+X_{2})=\dis\frac{3-q}{5-3q}(\sigma_{1}^{2}+\sigma_{2}^{2}); 1 \leqslant q<\frac{5}{3}$$
and $$\overline{V}_{2q-1}(X_{1}+X_{2})=\dis\frac{3-q}{q+1}(\sigma_{1}^{2}+\sigma_{2}^{2})=\overline{V}_{2q-1}(X_{1})+\overline{V}_{2q-1}(X_{2}); $$
In this case $cov(X_{1},X_{2})=0,$ because  $V(X_{1}+X_{2})=V(X_{1})+V(X_{2})+2cov(X_{1},X_{2}).$

\begin{cor}
Let $X_{1}, X_{2},...,X_{n}$ be n q-independent random variables with same q-Gaussian distribution $N_{q}(m,\sigma^{2}),$ then $\overline{X}_{n}=\frac{1}{n}\dis\sum_{i=1}^{n} X_{i}$ follows the q-Gaussian $N_{q}(m,\frac{\sigma^{2}}{n}).$	
\end{cor}

Observe that $V(\overline{X})  \longrightarrow 0 $ as $n \longrightarrow \infty$. Since $E(\overline{X})=m$, then the estimates of m becomes increasingly concentrated around the true population parameter. Such an estimate is said to be consistent.
\\The empirical variance
${S_{n}}^{2}=\frac{1}{n} \sum_{i=1}^{n}  (X_{i}-\overline{X_{n}})^{2}$ is not an unbiased estimate of $\sigma^{2}$.
Indeed
\begin{align*}
&\begin{aligned}{}
E({S_{n}}^{2})&=\dis\frac{n-1}{n} V(X_{1})\\
&=\dis\frac{n-1}{n} \dis\frac{3-q}{5-3q}\sigma^{2}\\
\end{aligned}&&
\end{align*} 
Therefore 
$$\widehat{\sigma}^{2}=\dis\frac{n}{n-1}\dis\frac{5-3q}{3-q}S_{n}^{2}.$$
is an unbiased estimate of ${\sigma}^{2}.$

\begin{prop} Law of Large Numbers (LLN):
If the distribution of the i.i.d. q-independent $X_{1}, . . . , X_{n}$ is such that $X_{1}$ has finite q-expectation,	i.e. $|E_{q}(X_{1})| < \infty$, then the sample average
$$\overline{X_{n}}= \frac{X_{1} +...+X_{n}}{n}  \longrightarrow E_{q}(X_{1})$$
converges to its expectation in probability.
\end{prop}

\begin{thm} Central Limit Theorem (CLT):\cite{4,2,8,5}
	For $q \in (1, 2),$ if $X_{1}, X_{2},.. X_{n}$ are q-independent and identically distributed with q-mean $m_{q}$ and a finite second $(2q- 1)$-moment $\sigma^{2}_{2q-1},$ then	
	$$Z_{n} = \frac{X_{1}+...+X_{n}-n m_{q}}{C_{q,n,\sigma}}$$  
	$q_{-1}$ converges to  $N_{q_{-1}}(0,1)$ Gaussian distribution.
\end{thm} 
Let $X$ be an arbitrary random variable with known variance $\dis\frac{3-q}{5-3q}\sigma^{2}$, and let $ X_{1}, X_{2}, . . . , X_{n}$ be n- q-independent  random variables with commun Gaussian distribution $N_{q}(m,\sigma^{2})$.

According to central limit theorem, the confidence interval for $m$  with level  $1-\alpha$ for arbitrary data and known $\sigma^{2}.$ is defined as 
$$m \in {(}\overline{X}_{n}\pm \frac{z_{1-\frac{\alpha}{2} }\sigma} {C_{q,N} \sqrt{n}}\Big{)}.$$
\\Where $z_{1-\frac{\alpha}{2}}$ is the quantile for $N_{q}(0,1)$ with level $1-\frac{\alpha}{2}:$
$$\int_{-\alpha}^{\alpha} N_{q}(0,1)(x)dx=1-\alpha,  \alpha \in ]0,1[.$$

\vskip0.2cm
\bibliographystyle{amsplain}

\end{document}